\documentstyle[12pt,amssymb,pb-diagram]{article}
\begin{document}
 
\title{On endomorphisms of projective bundles}
\author{Ekaterina Amerik
\thanks{Universit\'e Paris-Sud, Laboratoire des Math\'ematiques,
B\^atiment 425, 91405 Orsay, France. Ekaterina.Amerik@math.u-psud.fr} }

\date{September 13, 2002}
\maketitle

\def \C{{\bf C}}
\def \R{{\bf R}}
\def \Z{{\bf Z}}
\def \p{{\bf P}}
\def \N{{\bf N}}
\def \Q{{\bf Q}} 

\begin{abstract}

Let $X$ be a projective bundle. We prove that $X$ admits an endomorphism 
of degree
$>1$ and commuting with the projection to the base, if and only if $X$ 
trivializes after a finite covering. When $X$ is the 
projectivization of a vector
bundle $E$ of rank 2, we prove that it has an endomorphism of degree $>1$
on a general fiber only if $E$ splits after a finite base change.

\end{abstract}

It is clear that, for a complex projective variety $X$, the existence of 
endomorphisms $f:X\rightarrow X$ of degree bigger than one imposes
very strong restrictions on the geometry of $X$. On the other hand,
for any variety $B$, the product $B\times \p^r$ has such endomorphisms.
So it might be interesting to study the following question:

\
 
{\it Let $X$ be a projective bundle over a smooth projective complex
variety $B$, 
$p:X\rightarrow B$ the projection map.
When does $X$ admit a surjective endomorphism of degree bigger
  than one?}

\

This question is the subject of the present article.

Remark that any surjective endomorphism $f$ of $X$ is finite. Indeed, the
inverse image map $f^*$ on the rational cohomologies of $X$ is injective,
because $$f_{\ast}f^*=\deg(f)\cdot id,$$ so an isomorphism, but if $f$
contracts a curve, this map obviously cannot be surjective.

We will prove the following

\

{\bf Theorem 1} {\it X admits an endomorphism over $B$ (i.e. commuting
  with p) 
and 
of degree bigger than one if and only if $X$ trivializes after a
finite base change.}

\

In other words, the existence of such endomorphisms is equivalent
to the existence of a finite covering $B'$ of $B$ with covering
group
$G$, such that $X$ is a quotient of the product $B'\times \p^r$ by a 
suitable action of $G$ ( note that if $B''\rightarrow B$ is a
finite morphism,
then there exist a normal covering $B'\rightarrow B$ which factors through
$B''$, that is, the normalization of $B$ in the splitting field of the
extension $\C(B)\subset \C(B'')$ ).

It will be clear from the proof ( of the ``only if'' part ) 
that the covering $B'$ can be chosen
unramified, in other words, our theorem can also be reformulated as follows:

\

{\it X admits an endomorphism over $B$ of degree bigger than one 
if and only if $X$
is flat, coming from a representation 
$$\rho: \pi_1(B)\rightarrow PGL(r+1)$$ with finite image.}

\

For example, if $B$ is of general type, then some power of any endomorphism
$f$
of $X$ should be over $B$ and thus our theorem answers the question.
Indeed, $f$ sends fibers to fibers as $B$ is not covered by rational 
curves. Therefore $f$ must be over some endomorphism $\sigma$ of $B$,
i.e. $$ p\cdot f= \sigma \cdot p,$$
but $\sigma$ is of finite
order
and that means that for some $n$, $f^n$ is over $B$.

\

As for endomorphisms which are not over the base $B$, we have only
partial results for $\p^1$-bundles which we describe below.

Let $X=\p(E)$ with $E$ a vector
bundle of rank 2. It is clear that if $h$ is an endomorphism of 
$X=\p(E)$, then some power of it, $f=h^k$, is over an 
endomorphism $g$ of $B$ (i.e. $p\cdot f=g\cdot p$.) Indeed, one can 
characterize the fibers of $X$ as connected curves having zero 
intersection number with all $p^*L, L\in NS(B)$ (by the projection 
formula). That is, to check that $f=h^k$ sends fibers to 
fibers, it suffices (again by the projection formula) to show 
that $f^*(p^*NS(B)\otimes \Q)=p^*NS(B)\otimes \Q$. 
As $p^*NS(B)\otimes \Q$ is a hyperplane in $NS(X)\otimes\Q$ such that
any divisor class $D$
from this hyperplane satisfies $D^{\dim X}=0$, $h^*p^*NS(B)\otimes\Q$ 
is a hyperplane with the same property. But obviously there is only a 
finite number of such hyperplanes; $h$ interchanges them and so some 
power of $h$ must transform each such hyperplane into itself.

So we may suppose $X$ has an endomorphism $f$ of degree $>1$ sending 
fibers to fibers, and there are
two different cases:

1) $\deg(f|_{X_b})>1$ for a general fiber $X_b, b\in B $ (and hence, 
as it is easy to see, for all of them);

2) $\deg(f|_{X_b})=1$ for all $b\in B$.

The second case means of course that $g^*E\cong E\otimes L$ where $L$ 
is a line bundle. A simplest
example is that of $B$ an elliptic curve, $g$ the multiplication by an
integer and $E$ the only non-trivial extension of the structure sheaf
by itself: we have that $g^*E=E$ and $g$ extends to an endomorphism of $X$.
This case will not be further considered in this paper.

For the first case, we will prove 

\

{\bf Theorem 2:} {\it Suppose that $X=\p(E)$ with $E$ a vector bundle
of rank 2. If there exists an endomorphism of X which
is not of degree one on a general fiber, then $E$
splits into a direct sum of two line bundles after a finite
base change.}

\

The converse cannot of course be true as it is clear from the
discussion
above. The other disadvantage of this theorem is that we are not able
to say much about the finite base change which appears. One can nevertheless
obtain some more precise results by checking what does the proof of our
Theorem 2 give for some particular varieties.

If for example $B=\p^n$, we obtain the following

\

{\bf Proposition 3} {\it Suppose $B=\p^n$ and $X=\p(E)$ with $E$ a vector bundle
of rank 2. There exists an endomorphism of X of degree bigger than one 
if and only if $E$ is a direct sum of two line bundles.}
 
\

The reason is that in the proof of the Theorem 2 we arrive at the following 
statement (for any $B$): either $X$ trivializes after a finite unramified
covering of $B$, or $E$ has a subbundle. In the case $B=\p^n$, it means
of course that $E$ splits. The Proposition can also be obtained directly;
nevertheless, using the Theorem 2 
we see that the existence of endomorphisms implies the splitting of the
bundle for the whole class of bases $B$ which are ``similar to $\p^n$'', 
for example,
simply connected and such that $H^1(B,L)=0$ for any line bundle $L$.

Another remark on this Proposition is that projectivizations of split
bundles over $\p^n$ are exactly the projective bundles over $\p^n$
which are toric ( by a lemma of Druel, \cite{D} ). Toric varieties have
obvious endomorphisms. One can ask the following question: is it true that
a rational variety with an endomorphism of degree bigger than one, is toric?
Though I could not completely check it even for surfaces, 
I don't know any counterexample.

Throughout the paper, $\p(E)$ means the projective bundle of lines in $E$
(not hyperplanes in $E$), so that  $p_*{\cal O}_{\p(E)}(1)=E^*$.

I am very much indebted to L. Manivel and M. Rovinsky 
for several suggestions without which
the proof of Theorem 1 would be longer and less natural. 

\

{\bf 1. Proof of Theorem 1}

\

The ``if'' direction is as follows: 
if $X$ is a quotient of the trivial bundle by a
finite group $G$
as above, then to find an endomorphism of $X$ is the same as to find
an
endomorphism of $\p^r=\p(V)$ which is equivariant with respect to $G$.
$G$ acts on $S^m(V^*)$
 in the natural way. It is clear that for a big $m$ one
can find a $G$-invariant $f\in S^m(V^*)$ such that its zero locus
$D\in \p(S^m(V^*))$ is smooth (use Bertini theorem). 
Choose coordinates
$x=(x_0:x_1:\dots :x_r)$ and
consider the map $I:\p^r\rightarrow \p^r$:
$$I(x_0:x_1:\dots :x_r)=\left(\frac{\partial f}{\partial x_0}:
\frac{\partial f}{\partial x_1}:\dots:
\frac{\partial f}{\partial x_r}\right)$$ 
This is well-defined because $D$ is smooth, and in fact it can be
seen as
a $G$-equivariant map from $\p^r$ to its dual: differentiating
$f(gx)=f(x)$ (or, equivalently, referring to the Koszul complex), we see that the diagram below 
is commutative:

$$\begin{diagram} 
\node{\p^r} \arrow{e,t}{f} \arrow{s,l}{g} \node{\p^r} 
\arrow{s,r}{(g^{-1})^t} \\ 
\node{\p^r} \arrow{e,t}{I} \node{\p^r} \end{diagram}$$

Repeating this procedure for the dual action of $G$, we get a 
morphism \ 
$I':\p^r\rightarrow \p^r$ with $gI'(x)=I'((g^{-1})^tx)$, 
and the composition $I'I$ is the equivariant endomorphism 
which we are looking for.

\

Before proving the ``only if'' part of the Theorem,
let us state the following

\

{\bf Proposition 1.1} {\it Consider the parameter space $R^m(\p^r,\p^r)$ 
of regular
endomorphisms of $\p^r=\p(V)$ given by degree m polynomials, together with
the
natural action of $PGL(r+1)$:
$$g\cdot f(x)=gf(g^{-1}x).$$ If $m>r+1$, there exists an affine geometric 
quotient $R^m(\p^r,\p^r)/PGL(r+1)$.}

\

The proof is a sequence of lemmas:

\

{\bf Lemma 1.2} {\it $R^m(\p^r,\p^r)$ is affine.}

\

{\it Proof:} Indeed, $R^m(\p^r,\p^r)\subset \tilde{R}^m(\p^r,\p^r)$
where $\tilde{R}^m(\p^r,\p^r)=\p(V^*\otimes S^mV)$ is the space of 
rational self-maps of $\p^r$. The complement of  
$R^m(\p^r,\p^r)$ in $\tilde{R}^m(\p^r,\p^r)$ is the image of the incidence
variety $I\subset \p^r \times \tilde{R}^m(\p^r,\p^r)$ consisting of pairs
$(x,f)$ such that $f$ is not defined at $x$, by the projection $p_2$.
It is clear that $I$ is irreducible, $p_2$ is generically finite on $I$,
and the dimension count gives that $p_2(I)$ is a hypersurface.

\

{\bf Lemma 1.3} {\it If $m>r+1$, a unipotent element $U$ of $GL(V)$ does not
  stabilize
any element of $R^m(\p^r,\p^r)$.}

\

{\it Proof:}
To simplify the notations, we treat the case when $U$ is a single Jordan cell.
If $f$ is fixed by $U$, then by definition

$$
\begin{array}{c}
f_r(x_0, x_1,\dots,x_r)=\mu f_r(x_0+x_1, x_1+x_2,\dots,x_r)\\
f_{r-1}(x_0, x_1,\dots,x_r)+f_r(x_0, x_1,\dots,x_r)=\mu 
f_{r-1}(x_0+x_1, x_1+x_2,\dots,x_r)\\
\ldots \ldots\\
f_0(x_0, x_1,\dots,x_r )+f_1(x_0, x_1,\dots,x_r)
=\mu f_0(x_0+x_1, x_1+x_2,\dots,x_r)
\end{array}
$$
for some complex number $\mu$ ( where $\mu$ does not depend on $x$ 
because it can be seen as a regular function on $\p^r$, i.e. constant ).
Restrict to the projective line $x_2=x_3=\dots=x_r=0$: we get

$$g_r(x_0, x_1)=f_r(x_0, x_1,0,\dots,0)=\mu g_r(x_0+x_1, x_1),$$
so $g_r$ induces a constant function on the affine part $\{x_1\neq 0\}$
and so is of the form $a_{0,r} x_1^m$. Using induction on $i$, we get
that this is also true for 
$$g_{r-i}(x_0, x_1)=\frac{\partial^i}{\partial x_0^i}
f_{r-i}(x_0,x_1,0,\dots,0), 0\leq i\leq r.$$
So we have
$$
\begin{array}{c}
f_r(x_0, x_1,0,\dots,0)=a_{0,r} x_1^m\\
f_{r-1}(x_0, x_1,0,\dots,0)=a_{0,r-1} x_1^m+a_{1,r-1}x_1^{m-1}x_0\\
\ldots \ldots \\
f_{r-i}(x_0, x_1,0,\dots,0)=
a_{0,r-i} x_1^m+\dots+a_{i,r-i}x_1^{m-i}x_0^i\\
\ldots \ldots \\
f_0(x_0, x_1,0,\dots,0)=a_{0,0}x_1^m+\dots+a_{r,0}x_1^{m-r}x_0^r 
\end{array}
$$

for certain $a_{i,j}\in \C$.

So, if $m$ is bigger than $r$, then 
$f_i(1,0,\dots,0)=0$ for any $i$, a contradiction.

In the case when $U$ is not a single Jordan cell, the argument remains
almost the same, but we have to restrict on the projective
line $x_0=x_1=\dots=x_{i-1}=x_{i+2}=x_{i+3}=\dots=x_r=0$, 
where $i$ is such that the first non-trivial Jordan cell 
starts on the $i$th line.

\

{\bf Lemma 1.4} {\it The only one-parameter subgroup
$diag(\lambda^{a_0},\dots,\lambda^{a_r})\subset GL(V)$, $a_i\in \Z$,
which fixes
any element of $R^m(\p^r,\p^r)$, $m>1$, is the group of scalar matrices.}

\

{\it Proof:} Let $f$ be fixed by $diag(\lambda^{a_0},\dots,\lambda^{a_r})$. 
Without the loss of generality, we may assume $a_0=\min(a_i)$, 
$a_1=\max(a_i)$. Let $i_0$ be such that $f_{i_0}(1,0,\dots,0)\neq 0$ 
and analogously for $i_1$. We have
$$\lambda^{a_{i_0}}f_{i_0}(x_0,0,\dots,0)=\mu(\lambda)f_{i_0}
(\lambda^{a_0}x_0,0,\dots,0)=
\mu(\lambda)\lambda^{ma_0}f_{i_0}(x_0,0,\dots,0),$$
with $\mu(\lambda)$ a complex-valued function of $\lambda$,
from where $$\mu(\lambda)=\lambda^{a_{i_0}-ma_0}.$$ But in the same way
$$\mu(\lambda)=\lambda^{a_{i_1}-ma_1},$$ from where
$m(a_1-a_0)=a_{i_1}-a_{i_0}$, so, either $m=1$ or all $a_i$ are equal.

\

Now there are the classical results ( by Hilbert, Nagata, Mumford and others,
see for example \cite{PV} ) that if $U=Spec(A)$ 
is an affine variety and $G$ is a reductive
group acting on $U$, then the ring of invariants $Spec(A^G)$ is finitely
generated, and if moreover $G$ acts in such a way that all orbits are 
closed in $U$, 
then
$Spec(A^G)$ is a geometric quotient (that is, parametrizes the orbits).
 
By our lemmas, we are in this situation. 
Indeed, for any $f\in R^m(\p^r,\p^r),$
$Stab_{GL}(f)\subset GL(V)$
is an algebraic subgroup  which by Lemma 1.3 consists of
semisimple
elements. Take any of these elements and consider the minimal algebraic
subgroup of the stabilizer which contains it; 
this subgroup is an algebraic quasi-torus 
and so its connected component of the unity is diagonalizable.
 By Lemma 1.4, the only one-parameter subgroup 
in this group can be the center of $GL(V)$, i.e.  $Stab_{PGL}(f)\subset
PGL(V)$
is finite. So all the $PGL$-orbits are of the same dimension, i.e. closed, 
and so the Proposition 1.1 follows.

\

{\it Proof of the ``only if'' part of the Theorem:}

\

We may of course suppose that an endomorphism $f$ 
of $X$ over $B$ is given by the polynomials of degree $m>r+1$ on the fibers.
Consider the projective bundle $\tilde{R}^m_B(X,X)$ (the ``bundle of
rational
self-maps of $X$'') with fiber $\tilde{R}^m(\p^r,\p^r)$ and the
 transition functions $r_{\alpha, \beta}$ induced from the transition 
functions $g_{\alpha, \beta}$ of the projective bundle $X$ over $B$ 
by the rule 
$$r_{\alpha,\beta}f(x)=g_{\alpha, \beta}f(g_{\alpha,\beta}^{-1}x).$$
Consider the ``locus of regular maps'' 
$R^m_B(X,X)\subset \tilde{R}_B^m(X,X)$. 
There is a natural surjection $$\phi: R^m_B(X,X)\rightarrow
(R^m(\p^r,\p^r)/PGL(V))\times B$$ ( the ``quotient map''), 
which on the fibers is defined by choosing
an isomorphism $i_b:X_b\cong \p(V)$; this definition is independent
of $i_b$ as we factorized by all possible coordinate changes. 

The endomorphism $f$ is a section of $\tilde{R}^m_B(X,X)$ contained in 
$R^m_B(X,X)$. Its image $\phi(f)$ is a constant section $T\times B$
for some $T\in R^m(\p^r,\p^r)/PGL(V)$, because $R^m(\p^r,\p^r)/PGL(V)$
is affine. Let $t$ be a representative of $T$ in $R^m(\p^r,\p^r)$.

In a suitable trivialization $\{(U_\alpha, \psi_\alpha )\}$ of the bundle
$X$ over $B$,
we have
$f_\alpha=h_\alpha \cdot t$ for  
$h_\alpha$ a function on  $U_\alpha$ with values in $PGL(r+1)$.
Here ``$\cdot$'' means the natural action of $PGL(r+1)$
on $R^m(\p^r,\p^r)$: $(g\cdot f)(x)=gf(g^{-1}x).$

Refining $\{U_\alpha\}$ if necessary (in particular choosing them
simply 
connected), 
we can make 
$h_\alpha$ holomorphic. In fact, consider the subset 
$H\subset U_\alpha
\times PGL(r+1):$ 
$$H=\{(x,h): f_{\alpha ,x}=h\cdot t\}$$
which is obviously a subvariety. 
The map $u: H\rightarrow U_\alpha$ is finite,
and it is clear that one can construct holomorphic sections $h_\alpha$ over
a suitable refinement of $\{U_\alpha\}$ provided that $u$ is smooth.
The last assertion is true because $u$ is obtained from a smooth map
$PGL\rightarrow PGL(r+1)\cdot t$ by the base change 
$U_\alpha\rightarrow PGL(r+1)\cdot t$ given by $f_\alpha$.

So for $U_\alpha$, $U_\beta$ from our covering, 
$$\begin{array}{c} 
f_\alpha=h_\alpha \cdot t, h_\alpha\in PGL_{r+1}({\cal O}(U_{\alpha}))\\
f_\beta=h_\beta \cdot t, h_\beta\in PGL_{r+1}({\cal O}(U_{\beta})),
\end{array}$$
and on the intersection $U_{\alpha \beta}$, $f_\alpha=g_{\alpha \beta}
\cdot f_\beta$ where $g_{\alpha \beta}$ are just the transition 
functions for $X$. We therefore have

$$h_\alpha^{-1}g_{\alpha \beta}h_\beta\in Stab_{PGL}(t),$$
so by changing the trivialization we can make $g_{\alpha \beta}$ constant
functions with values in a finite group $Stab_{PGL}(t)$, q.e.d.

\

{\bf Remark:} The case of ruled surfaces over a curve is particularly 
nice because the finite subgroups of $PU(2)$ are well known:
cyclic, dihedral, $A_4$, $S_4$, $S_5$. For each such subgroup $G$ the number
of corresponding ruled surfaces can be computed starting with the 
following classical formula (see e.g. \cite{S}, 7.2, exercises):
$$\# Hom(\pi_1(C), G)=\#G\sum_{r\in R}\left(\frac{\#G}{rk(r)}\right)^{2g-2}$$
where $g$ is the genus of $C$ and $R$ 
is the set of irreducible representations 
of $G$.

\

{\bf 2. Proof of Theorem 2}

\

Following the remarks from the Introduction, 
we assume that our endomorphism $f:X\rightarrow X$ satisfies
$p\cdot f=g\cdot p$ for some $g\in End(B)$. 

\

The following lemma is well-known, see for example \cite{H}, II.6:

\

{\bf Lemma 2.1} {\it The following problems are equivalent:

1) To find a morphism $f:X\rightarrow X$, $p\cdot f=g\cdot p$;

2) To find a surjective bundle
map
$$p^*g^*E^*\rightarrow L$$ with $L$ a line bundle on $X$; 
the line bundle $L$ is then $f^*({\cal O}_{\p(E)}(1))$; 

3) To find a morphism 
$$\phi:X=\p(E)\rightarrow \p(g^*E)$$ over $B$, 
$\phi^*({\cal O}_{\p(g^*E)}(1))=L.$}

\

{\bf Lemma 2.2} {\it The ramification locus of $\phi$ as above
cannot contain fibers of $X$.}

\

{\it Proof:} This is just a statement that if a morphism
$\p^1\rightarrow\p^1$
is given by polynomials of degree $k$, then its degree is exactly $k$
and not less, which is obvious.

\

We will prove a statement which is slightly more general than the Theorem 2:

\

{\bf Theorem 2A} {\it Let $E$, $F$ be vector bundles of rank 2 on $B$, 
$X=\p(E)$, $Y=\p(F)$. If there exist a morphism 
$\phi: \p(E)\rightarrow \p(F)$
over $B$ which is of degree bigger than one, then $E$ and $F$ both split
after a finite base change.}

\

We will follow the strategy of the proof of Theorem 1.
First, let $V,W$ be vector spaces of dimension two. 
Consider the affine variety $R^m(\p(V),\p(W))$ of morphisms from $\p(V)$
to $\p(W)$ given by degree $m$ polynomials, together with the action
of $PGL(V)\times PGL(W)$:
$$(g,h)\cdot f(x)=hf(g^{-1}x).$$
In this case, the quotient (i.e. the spectrum of the ring of invariants)
is not a geometric quotient: the orbits are
not all closed and so are not separated by the invariants.
However, it is easy to understand when different orbits
correspond to the same point of the quotient, by the following Proposition:

\

{\bf Proposition 2.3} 

{\it 1) A point $f\in R^m(\p(V),\p(W))$ has infinite
stabilizer in $PGL(V)\times PGL(W)$, if and only if in suitable
coordinates, $f$ is written as $y_0=x_0^m, y_1=x_1^m$.

2) An orbit which contains the map $y_0=x_0^m, y_1=x_1^m$ in its closure,
consists of maps which in suitable coordinates are of the form
$y_0=f_0(x_0, x_1), y_1=x_1^m$.}

\

{\it Proof:} 1) One can make the same calculations as in Proposition 1.1 : first
check that a unipotent $(g,h)\in GL(V)\times GL(W)$ does not stabilise any $f\in R^m(\p(V),\p(W))$,
supposing $g$ and $h$ in Jordan form and making a calculation (of the same
kind as before). Then we look for morphisms stabilized by an
one-parameter subgroup $(diag(\lambda^{c_0},\lambda^{c_1}),
diag(\lambda^{b_0},\lambda^{b_1}))$. Let $f$ be such a morphism given by
$$
\begin{array}{c}
y_0=f_0(x_0,x_1)=\Sigma_{i=0}^{m} a_ix_0^ix_1^{m-i}\\
y_1=f_1(x_0,x_1)=\Sigma_{i=0}^{m} a'_ix_0^ix_1^{m-i}

\end{array}$$

From the invariance of $f$ by the subgroup, we get 
$b_0-c_0i-c_1(m-i)$ are equal for all $i$ such that $a_i\neq 0$, 
and they are moreover equal to
$b_1-c_0i-c_1(m-i)$ for all $i$ such that $a'_i\neq 0$.
 
If $c_0=c_1$, then also $b_0=b_1$ and this is trivial in 
$PGL(V)\times PGL(W)$. If $c_0\neq c_1$, then at most one of the 
$a_i$'s and one of the $a'_i$'s is nonzero. As our $f$ is a morphism 
of degree $m$, we get, up to rescaling and permutation of $x_0$ and
$x_1$, $f_0(x_0,x_1)=x_0^m$ and $f_1(x_0,x_1)=x_1^m$.

Alternatively and more geometrically, one can remark that if a map $f$ 
is fixed by $(g,h)\in PGL(V)\times PGL(W)$, then $h$ preserves its 
ramification locus and $g$ preserves its branch locus. After checking 
that no positive-dimensional part of the stabilizer can be of the form 
$G\times {id}$ or ${id}\times G$, we conclude that both loci consist
of two points, which in suitable coordinates are $(0:1)$ and $(1:0)$. 
The ramification multiplicity at the two ramification points must be 
of course $m-1$ (as it cannot be bigger); we have set-theoretically
$f^{-1}((0:1))=(0:1)$ and the same for $(1:0)$. Which is what we had to show.

\

2) Suppose that $h$ is in a non-closed orbit.
By a generalization of the Hilbert--Mumford criterion 
(see \cite{B}, Theorem 4.2, or \cite{PV}, Theorem 6.9), 
we can reach its border by a suitable one-parameter subgroup 
$\{G_{\lambda}, \lambda \in \C^*\}$. We may consider 
$G_{\lambda}\in SL(V)\times SL(W)$. Choose coordinates such that 
all $G_{\lambda}$ diagonalize:

$$G_{\lambda}=(diag (\lambda^{-c},\lambda^{c}),\quad 
diag(\lambda^{-b},\lambda^{b})).$$

If $h$ is given by 
$$
\begin{array}{c}
y_0=h_0(x_0,x_1)=\Sigma_{i=0}^{m} a_ix_0^ix_1^{m-i}\\
y_1=h_1(x_0,x_1)=\Sigma_{i=0}^{m} a'_ix_0^ix_1^{m-i}
\end{array},$$

we have $G_{\lambda}\cdot h$ as follows:

$$
\begin{array}{c}
y_0=\Sigma_{i=0}^{m} a_i \lambda ^{(2i-m)c-b}x_0^ix_1^{m-i}\\
y_1=\Sigma_{i=0}^{m} a'_i\lambda ^{(2i-m)c+b}x_0^ix_1^{m-i}
\end{array}
$$

If $c=0$ but $b\neq 0$, then in the limit of $G_{\lambda}\cdot h$ as 
$\lambda$ goes to zero, we get a constant morphism (from $\p^1$ to $0$ or
$\infty$, depending on the sign of $b$ (so not an element of $R^m$). Suppose that $c>0$ (in the case $c<0$ it will be the same computation),
so that
the exponents $(2i-m)c-b$ and $(2i-m)c+b$ grow together with $i$.
Consider $i_1=\min\{i, a'_i\neq 0\}$, and let $K=(2i_1-m)c+b$.
 It is clear that if we want to get a non-constant map in the limit, then
we must have, for $i_0=\min\{i, a_i\neq 0\}$, that $(2i_0-m)c-b=K$, and
the map in the limit will be 
$$
\begin{array}{c}
y_0=a_{i_0}x_0^{i_0}x_1^{m-i_0}\\
y_1=a'_{i_1}x_0^{i_1}x_1^{m-i_1}
\end{array}
$$

and, as we want to stay in $R^m$, either $i_0=0, i_1=m$, or vice versa.

\

{\it Proof of Theorem 2A}

\

Let $M$ be the (affine) quotient $R^m(\p(V),\p(W))/PGL(V)\times PGL(W)$,
and let $M^0$ be its subset parametrizing closed orbits. Let $X=\p(E)$,
$Y=\p(F)$, and consider the bundle $R^m(X,Y)$ and the surjection
$\phi: R^m(X,Y)\rightarrow M\times B$ as in the proof of Theorem 1.
The image $\phi(f)$ is a constant section $T\times B$. If $T\in M^0$,
then we repeat the argument of Theorem 1 to conclude that both $X$ and $Y$
trivialize after a finite unramified covering. If $T\in M-M^0$, then the 
restriction 
of $f$ to each fiber has a ramification point of multiplicity $m-1$, so the
ramification divisor of $f$ has a component of multiplicity $m-1$.
There are two cases possible: 

In the first case, this component is the only
component of the ramification. Then the restriction of $f$ on each fiber 
is ramified at two points only, with multiplicity $m-1$ at each (as obviously
a morphism from $\p^1$ to $\p^1$ cannot have only one ramification point).
The ramification is a double covering of the base. After pulling back our
projective bundle $\p(E)$ to the ramification, we have two non-intersecting
sections, giving the splitting of the pull-back.

In the second case, the ramification has an irreducible component $R_1$ of 
multiplicity $m-1$ which is a section of $\p(E)$. We claim it does not 
intersect the rest $R_2$ of the ramification divisor.
Indeed, if it does, then any intersection point with $R_2$ will
be a ramification point of multiplicity at least $m$ of a morphism
from $\p^1$ to $\p^1$ of degree $m$, which is impossible. Let $B'$ be an
irreducible component of $R_2$. The pullback of $X$ to $B'$ has two
non-intersecting sections, so is the projectivization of a split vector
bundle.

The bundle $F$ also splits after a finite base change. Indeed, in the 
first case the branch locus is a double covering of the base, and in the
second case the image of $R_1$ is a section and the image of $R_2$ is a 
multisection which does
not intersect the image of $R_1$.

\

Theorem 2 follows, of course. Unfortunately,
the proof does not say much about the base change 
$B'\rightarrow B$ after which
$E$ would split; it is not even clear that one can choose $B'$ non-singular.
However the situation is better in some particular cases:

\

{\bf Proposition 2.4} 
{\it If $E$ is semistable (in the usual sense if $B$ is a curve, or 
$H$-semistable with respect to some
 ample divisor $H$ if $\dim(B)>1$), 
then the base change $B'$ can be 
chosen smooth and unramified over $B$.}

\

{\it Proof:}

Let $\phi$ be a morphism from $X=\p(E)$ to $\p(F)$ over $B$, and let
$\xi$ resp. $\xi'$ be the tautological divisor classes on $X$ resp. $Y$. Let $k$ be
the
degree of $\phi$ (on the fibers). By Lemma 2.1, this morphism gives
a non-vanishing section of $p^*F(L)$ on $X$,  where $L=k\xi+p^*D$
is the inverse image of $\xi'$. So in the Chow ring of $X$ we must have
$$c_2(p^*F(L))=c_2(p^*(F))+k\xi c_1(p^*F)+p^*D c_1(p^*F)+k^2\xi^2 +
2k\xi p^*D + p^*D^2=0.$$

Substituting $\xi^2=-\xi c_1(p^*E)-c_2(p^*E)$, and taking into account
the absence of linear relations between the elements of $p^*CH^2(B)$ and 
$\xi p^*CH^1(B)$, we get

$$D\equiv \frac{1}{2}(k\cdot\det E-\det F)$$
modulo numerical equivalence.

Now the canonical classes are

$$\begin{array}{l} K_{\p(E)}=-2\xi+p^*K_B-p^*(\det E)\quad\mbox{and}\\
K_{\p(F)}=-2\xi'+p^*K_B-p^*(\det F),\end{array}$$ 
so the ramification divisor is, up to numerical equivalence,
$$K_{\p(E)}-\phi^*K_{\p(F)}=(2k-2)\xi+(k-1)p^*\det E.$$

Remark that the determinant of $$S^{2k-2}E^*((k-1)\det E)=p_*
{\cal O}((2k-2)\xi+(k-1)p^*\det E)$$ is trivial.

Let $M\in |l\xi+p^*A|$ be a component of the 
ramification divisor. By direct image, we have 
$$H^0(X, {\cal O}(l\xi+p^*A))=H^0(B,S^l(E^*)\otimes {\cal O}(A)).$$ $E$ being 
$H$-semistable,
so is $S^l(E^*)\otimes {\cal O}(A))$ ( \cite{HL}, Theorem 3.1.4 ). 
Computing determinants, we have that if this bundle has a section, 
then $AH^{n-1}\geq \frac{l}{2}c_1(E)H^{n-1}$. As this is true for 
any component of the ramification, we must in fact have equality. 

We now can prove that $M$ is smooth and unramified over $B$. To check this
at the point $x\in M$, consider the restriction to a curve 
$C=H_1\bigcap H_2 \bigcap\dots\bigcap H_{n-1},$
where $H_i$ are sufficiently general divisors containing $x$ and 
from $|rH|$ with $r$ big enough
(we want $C$ to be a ``sufficiently general curve passing through $x$'' and
to be in the cohomology class proportional to $H^{n-1}$).
By standard Bertini-type statements, $C$ is smooth and the restriction
$M_C=M\bigcap p^{-1}(C)$ is a reduced (but possibly not irreducible) 
divisor in the ruled surface
$S=p^{-1}(C)$. If we show that $M_C$ is smooth and unramified over $C$, it
will follow that $M$ is smooth and unramified over $B$ at the point $x$.

The statement about $M_C$ is implied by the following 
simple lemma on ruled surfaces:

\

{\bf Lemma 2.5} {\it Let $S=\p(G)$ be a ruled surface over the base curve $C$,
and let $\xi_C$ be the universal line bundle.
Let $N$ be an irreducible divisor on $S$, which is numerically
$m\xi_C+\frac{m}{2}p^*\det G+p^*D$, $m\geq 1$, with $D$ a divisor on $C$.

\

1) If $\deg(D)<0$, then $m=1$ (i.e. $N$ is a section);

2) If $\deg(D)=0$, then $N$ is smooth and unramified over $C$.}

\

This lemma (in a somewhat different form) appears
in \cite{H} ( Chapter V, Proposition 2.21), so we will not prove it here 
(the argument is to compare the
arithmetic genus of $N$ obtained by adjunction formula, with its
geometric genus obtained by Hurwitz' formula).

\ 

Now let $G=E|_C$; we have that $M_C$ is numerically equivalent to
$l\xi_C+\frac{l}{2}p^*\det G$. If $M_C$ is irreducible, we are done.
Notice that in any case $M_C$ is connected (by a well-known connectedness
theorem, see for example \cite{FL}, Theorem 2.1 ). In particular $M_C$ cannot split
into irreducible components of the same numerical type
$m\xi_C+\frac{m}{2}p^*\det G$, as elementary calculation shows that they
don't intersect. So, if $M_C$ is not irreducible, it must have a component
of numerical type $\xi_C+\frac{1}{2}p^*\det G+p^*D$, $\deg(D)<0$. If we compute
its intersection number with the rest of $M_C$, we get $(l-2)\deg(D)$; this
is non-positive so it is a contradiction. Thus $M_C$ is irreducible,
smooth and unramified over $C$, q.e.d.

\

{\bf Remark} Trying to apply the similar considerations to the case
when $E$ is unstable with respect to all ample divisors, in order to study
$R_1$ and $R_2$ appearing in the end of proof of Theorem 2A, we can obtain
the followng: if the reduction of $R_2$ is not a section, then the
reduction of $R_1$ is a section which comes from the destabilizing subsheaf
of $E$ (so that the destabilizing subsheaf is a subbundle independent of
the polarisation). 

\

{\it Proof of  Proposition 3}

\

Let $B=\p^n$.
Let $X=\p(E)$ have endomorphisms of degree bigger than one. We must 
first treat the case not covered by Theorem 2A, that is, the case of 
endomorphisms being of degree one on fibers. If $X$ has such an 
endomorphism $h$, some power of $h$ is over an endomorphism $g$ of 
$\p^n$ and for this $g$ we have $g^*E=E\otimes L$, $L\in Pic(\p^n)$. 
We claim that it is possible only if $E$ is a shift of the trivial 
bundle. One way to see this is the following: let $k$ be such that 
$g^*h=kh$ for $h$ the hyperplane section divisor. The calculation 
of Chern classes gives that $L=\frac{1}{2}(k-1)\det E$,
and $c_1^2(E)=4c_2(E)$.
So by Bogomolov's inequality (see \cite{Bog}, or \cite{OSS}, end of 
II.2, for our particular case of bundles on $\p^n$) $E$ is not stable, 
and if it is semistable, then a destabilizing subsheaf is a subbundle, 
that is, $E$ is a shift of the trivial bundle. In the case when $E$ is 
unstable, let $F\subset E$ be the maximal destabilizing subsheaf; the 
maximal destabilizing subsheaf of $E\otimes L$ is then $F\otimes L$. 
But $F^{\otimes k}=g^*F$ also destabilizes $E\otimes L$, so 
$kc_1(F)\leq c_1(F)+c_1(L)$. From $c_1(F)> \frac{1}{2}c_1(E)$ and 
$L=\frac{1}{2} (k-1)\det E$, we obtain a contradiction if $k\neq 1$.

So we may assume that $f$ is not of degree one on a general 
fiber. From the proof of Theorem 2A we see that either $X$ becomes 
a trivial bundle after a finite unramified covering, which means 
that $X$ already is trivial, or $E$ has a subbundle. As on $\p^n$, a 
rank-two bundle with a subbundle splits, the ``only if'' direction follows.

The ``if'' part comes from the fact that the projectivization of
a direct sum of line bundles over $\p^n$ is toric, and toric varieties
have obvious endomorphisms coming from homotheties of the associated fans.

\

{\bf Remark:} For $V,W$ vector spaces of higher dimension, 
we have many non-closed orbits
under $PGL(V)\times PGL(W)$ in $R^m(\p(V),\p(W))$. What one does obtain,
is that the non-closed orbits consist of maps ``with an invariant linear
subspace'': there exist linear subspaces $U_1\subset \p(V)$, 
$U_2\subset \p(W)$,
such that $f^{-1}(U_2)=U_1$. One then might try to proceed by induction.
However, I could not work it out even in dimension 3.

{}

\end{document}